\documentclass{amsart}

\usepackage{amsmath,amsthm,amssymb,IMjournal}
\makeatletter 

\@addtoreset{equation}{section}
\makeatother 

\newtheorem{lemma}{Lemma}[section]

\newtheorem{remark}{Remark}[section]
\usepackage{color}
\allowdisplaybreaks 
\numberwithin{equation}{section}
\usepackage{cases}
\usepackage{mathabx}
\usepackage{epsf, graphicx}
 \usepackage{ulem}
\usepackage{bm}

\begin{document}
\newtheorem{The}{Theorem}[section]

\numberwithin{equation}{section}

\title{Explicit Finite Element Error Estimates for nonhomogeneous Neumann problems }

\author{\|Qin |Li|, Beijing, China,~~
        \|Xuefeng |Liu|, Niigata, Japan}

\rec { March 31, 2018}

\dedicatory{Cordially dedicated to ...}

\abstract
The paper develops an explicit a priori error estimate for finite element solution to {nonhomogeneous Neumann problems}.
For this purpose, the hypercircle equation over finite element spaces is constructed and
the explicit upper bound of the constant in the trace theorem is given.
Numerical examples are shown in the final section, which implies the proposed error {estimate} has the convergence {rate} as $0.5$.
\endabstract

\keywords
   Finite element methods; {nonhomogeneous Neumann problems}; explicit error estimates
\endkeywords

\subjclass
65N15,65N30
\endsubjclass

\thanks
   The research has been supported by the National Science Foundations of China (No. 11426039, 11571023, 11471329) and 
   Research Foundation for Youth Scholars of Beijing Technology and Business University (QNJJ2014-17) for the first author and
   Japan Society for the Promotion of Science, Grand-in-Aid for Young Scientist (B) 26800090 and 
   Grant-in-Aid for Scientific Research (C) 18K03411, for the second author.
\endthanks

\section{Introduction}\label{sec1}

The Steklov type differential equation problem involves the {Neumann} boundary conditions{. It models} various physical phenomenon{, for example, the
vibration modes of a structure in contact with an incompressible fluid \cite{BK}, the antiplane shearing on a system of collinear faults under slip-dependent friction law \cite{BII}}.
There is {wide} literature on numerical schemes to solve this type { of problems} by using for example, finite element method (FEM); see \cite{BrO, LLX}.
Also, the Steklov type eigenvalue problem is a fundamental problem in mathematics.
For example, the  {optimal} constant appearing in the trace theorem for Sobolev space{s} is {given by the smallest eigenvalue of a Steklov type} eigenvalue problem  {raised to the power $-\frac{1}{2}$; see e.g., \cite{SEVE}}.
Efforts have been {made} on bounding eigenvalues by using conforming or non-conforming FEMs; see \cite{LLX, YLL}.

Most of the existing literature focuses on the convergence analysis of discrete solution,
while there has been very rare work on the explicit bound of the solution error.
Recently, in the newly developed field of verified computing, the quantitative error {estimate} (e.g., explicit values of error) is
desired. For example, the explicit values or bounds of the error constants are required in solution verification of
non-linear partial differential equations; see, e.g., \cite{Takayasu}.

In this paper, we apply the finite element method to solve the Steklov type differential equation and
provide a priori error {estimate} for the FEM solution.
{
The main idea in developing a priori error estimaion can be regarded as a direct extension of the one proposed by Liu in \cite{LO}, where
the a priori error estimaion is constructed by using hypercircle equation for homogeneous boundary condtions.
Such ideas can be further tracked back to the one of Kikuchi in \cite{Kikuchi-Saito}, where a posteriori error estimation is considered.}
{This} a priori {estimate} can be used {for} bounding eigenvalue under the framework proposed by \cite{Liu} {and it will be the topic of a forthcoming paper}.

The rest of this paper is organized as follows. In section 2, we describe the problem to be considered.
In section 3, we construct the hypercircle equations over FEM spaces, based on which we deduce computable error estimates. 
In section 4, we discuss the constant appearing in the trace theorem and
propose the explicit {a} priori error {estimate} for {nonhomogeneous Neumann problems}.
In section 5, the computation results are presented.

\section{Preliminaries}\label{sec2}
Throughout this paper, we use the standard notation (see, e.g. \cite{BaO}) for the Sobolev spaces $H^m(\Omega)$ ($m>0$).
The Sobolev space $H^0(\Omega)$
coincides with $L^2(\Omega)$.
Denote by $\|v\|_{L^2}$ or $\|v\|_{0}$ the $L^2$ norm of $v\in L^2(\Omega)$;
$|v|_{H^m(\Omega)}$ and $\|v\|_{H^m(\Omega)}$ the seminorm and norm {in} $H^{m}(\Omega)$ , respectively.
{Symbol} $(\cdot, \cdot)$ {denotes} the inner product in $L^2(\Omega)$ or $(L^2(\Omega))^2$.
The space $H(\mbox{div}, \Omega)$ is defined by
$$
H(\mbox{div}, \Omega):=\{q\in (L^2(\Omega))^2\mid \text{div }\:q\in L^2(\Omega) \}.
$$

{ W}e are concerned with the following model problem
\begin{equation}\label{1}
\left\{
\begin{array}{rcl}
-\Delta u+ u &=& 0,\ \ \  {\rm in}\  \Omega\:, \\
\frac{\partial u}{\partial { \bm{n}}}&=&f,\ \ \ {\rm on}\ \Gamma=\partial\Omega \:,
\end{array}
\right.
\end{equation}
where $\Omega \subset \mathcal{R}^{2}$ is a bounded polygonal domain, $\frac{\partial }{\partial { \bm{n}}}$ is the outward normal derivative on boundary $\partial\Omega$.

A weak formulation of the above problem is to find $u \in
V=H^{1}(\Omega)$ such that
\begin{eqnarray}\label{11}
a(u,v)&=& b(f,v)~~~\quad\forall v\in V
\end{eqnarray}
where
$$
a(u,v)=\int_{\Omega}\big(\nabla u\nabla v+uv\big)dx,~~~ b(f,v)=\int_{\partial\Omega}fvds\:.
$$
Also, define $~\|u\|_{b}=b(u,u)^{1/2}$.


We also have the following regularity {result for the solution} of problem (\ref{1}); see, for example, \cite{GRI}.
\begin{lemma}
If $f\in L^{2}(\partial\Omega)$, then
$u\in H^{1+\frac{r}{2}}(\Omega)$;
if $f\in H^{\frac{1}{2}}(\partial\Omega)$, then $u\in
H^{1+r}(\Omega)$; here, $r\in (\frac{1}{2}, 1]$, especially $r=1$ when $\Omega$ is convex and $r<\frac{\pi}{\omega}$ (with $\omega$ being the
largest inner angle of $\Omega$) otherwise.
\end{lemma}

\paragraph{\bf Finite element approximation}

Let $\mathcal{T}_{h}$ be a {shape regular} triangulation of {the} domain $\Omega$.
For each element $K\in \mathcal{T}_{h}$, denote by $h_{K}$ the longest edge length of $K$ and define the mesh size $h$ by
$$h:=\max\limits_{K\in \mathcal{T}_{h}}h_{K}.$$
Define {by $E_{h}$} the set of edges of the triangulation and $E_{h,\Gamma}$ the {set of }edges on the boundary of {$\Omega$}.
{The finite element space $V^{h} (\subset V)$ consists of piecewise linear and continuous functions}. Assume that dim$(V^{h})=n$. The conforming finite element approximation of (\ref{11}) is defined as follows: Find $u_{h} \in V^{h}$ such that
\begin{eqnarray}\label{12}
a(u_{h}, v_{h})&=&b(f,v_{h})~~~\forall v_{h}\in V^{h}.
\end{eqnarray}

In this paper, the following classical finite element spaces will also be used in constructing the a priori {estimate}.

(i)  Piecewise constant function space{s} $X^{h}$ and $X_{\Gamma}^{h}$ {are defined as}:
\begin{eqnarray*}
&& X^{h}:=\{v\in L^2(\Omega)\mid v~\mbox{is constant on each element $K$ of }\mathcal{T}_{h} \}\\
&& X_{\Gamma}^{h}:=\{v\in L^2(\Gamma)\mid v~\mbox{is constant on each edge e} \in E_{h,\Gamma}\}.
\end{eqnarray*}

(ii) Raviart-Thomas FEM space {$W^{h}$}:
$$
W^{h}:=\{p_{h}\in H(\mbox{div }, \Omega)\mid p_{h}=(a_{K}+c_{K}x, b_{K}+c_{K}y)~\text{in}~ K\in \mathcal{T}_{h} \}
$$
where $a_{K}, b_{K}, c_{K}$ are constants on element $K$.

{The space} $W^{h}_{f_h}$ is a shift of $W^h$ corresponding to $f_{h}\in X^{h}_{\Gamma}$:
$$
W^{h}_{f_h}:=\{p_{h}\in W^{h} \: | \: p_{h}\cdot { \bm{n}}=f_{h}\in X_{\Gamma}^{h}~\text{on}~\Gamma \}
$$


\section{Hypercircle Equations}\label{sec3}
In this section, we first present two hypercircle equations which can be used to facilitate the error {estimate}.

Consider the boundary value problem: \footnote{The boundary condition can be extended to mixed one. For exmaple,
$\Omega=\Gamma_1\cup \Gamma_2$, ${\partial u}/{\partial { \bm{n}}}=f_1$ on $\Gamma_1$, $u=f_2$ on $\Gamma_{2}$.
}
\begin{equation}\label{100}
\left\{
\begin{array}{rcl}
-\Delta u+\alpha u &=& g\ \ \  {\rm in}\  \Omega \:, \\
\frac{\partial u}{\partial { \bm{n}}}&=&f\ \ \ {\rm on}\ \Gamma \:,
\end{array}
\right.
\end{equation}
with $\alpha$ being {a} positive constant and {$g\in L^2(\Omega)$}.
{A weak formulation of the above problem is to find $u \in
V=H^{1}(\Omega)$ such that
\begin{eqnarray}\label{1101}
\int_{\Omega}\big(\nabla u\nabla v+\alpha uv\big)dx&=& \int_{\Omega}gvdx+b(f,v)~~~\quad\forall v\in V
\end{eqnarray}
}
Corresponding to problem (\ref{100}), the following hypercircle equation holds; see, e.g., {page} 185 of \cite{Braess}.
\begin{The}
{Let $u$ be solution to problem (\ref{1101})}. For $v\in H^{1}(\Omega)$ and $v=0$ on $\Gamma$,
suppose that $\sigma \in H(\rm{div },\Omega)$ satisfies
$$
\sigma\cdot \vec{n}=f~on~\Gamma ~\text{ and }~
{\rm{div }}\sigma+g=\alpha v.
$$
Then, we have,
\begin{eqnarray}
\|\nabla(u-v)\|_{0}^2+\|\nabla u-\sigma\|_{0}^2+2 \alpha\|u-v\|_{0}^2=\|\nabla v-\sigma\|_{0}^2
\end{eqnarray}
\end{The}
\proof
The expansion of $\| \nabla v-\sigma\|^2=\|(\nabla v - \nabla u) + (\nabla u - \sigma)\|^2$ tells that
$$
\|\nabla v-\sigma\|^2 =\|\nabla v-\nabla u\|_{0}^2+\|\nabla u-\sigma\|_{0}^2+2( \nabla u-\sigma, \nabla (v-u))\:.
$$
Let $w:=v-u$.
From the definition of $u$ in (\ref{1101}), we have
\begin{equation}
\label{eq:hypercircle-a}
(\nabla u, \nabla w)= b(f,w) + \int_\Omega  (g - \alpha  u ) ~ w ~dx \:.
\end{equation}
Also, by applying Green's theorem to the term with $\sigma$, we have
\begin{equation}
\label{eq:hypercircle-b}
(\sigma, \nabla w)=\int_{\partial\Omega} (\sigma\cdot {n})~ w ds - \int_{\Omega} {\rm {div }}\sigma ~w ~ dx \\
= b(f,w) - \int_{\Omega}  ( \alpha v- g) ~w ~ dx\:.
\end{equation}
By taking (\ref{eq:hypercircle-a})-(\ref{eq:hypercircle-b}), we have $(\nabla u-\sigma, \nabla (v-u)) = \alpha \|v-u\|_0^2$, which
leads to the conclusion of this theorem.
\endproof

However, it is usually difficult to construct $\sigma$ { such that} ${\rm{div }}\sigma+g=\alpha v$ hold for general $v$ and $g$.
Below, we establish revised hypercircle equation over finite element spaces.
As a {preparation}, let us introduce two projection operators: $\pi_{h}$ and $\pi_{h,\Gamma}$.

\begin{itemize}

\item
For $g\in L^2(\Omega)$, define the projection $\pi_{h}$: $L^2(\Omega)\mapsto X^{h}$, such that,
$$(g-\pi_{h}g, v_{h})=0\quad\quad \forall v_{h}\in X^{h}.$$
The error {estimate} for $\pi_h$ is given by
\begin{eqnarray}\label{201}
 \|g-\pi_{h}g\|_{0}\leq C_{0}h|g|_{H^{1}(\Omega)}\quad{\forall g \in H^1(\Omega)} .
\end{eqnarray}
Here $C_0${:=$\max\limits_{K\in \mathcal{T}_{h}} C_0(K)/h$} depends on the triangulation and has {an} explicit upper bound.
For example, In {\cite{LS_2010,liu-kikuchi-tokyo-u}}, it is shown that the optimal constant is given by $C_0(K) = h_K/j_{1,1}$,
where $j_{1,1} \approx 3.83171 $ denotes the first positive root of the Bessel function $J_1$.
Upper bounds of $C_0$ for concrete triangles can be found in, e.g., \cite{Kobayashi, LiuFumio, liu-kikuchi-tokyo-u}.


\item
For $f\in L^2(\Gamma)$, define the projection $\pi_{h,\Gamma}$: $L^2(\Gamma)\mapsto X_{\Gamma}^{h}$,
$$
b(f-\pi_{h,\Gamma}f, v_{h})=0\quad\quad \forall v_{h}\in X_{\Gamma}^{h}.
$$


\end{itemize}


\begin{The}
Given $f_h \in X_\Gamma^{h}$, let $\tilde{u}\in
V$ and $\tilde{u}_{h} \in V^{h}$ be solutions {to} the following variational problems, respectively{:}
\begin{eqnarray}
&& a(\tilde{u},v)=b(f_{h},v)~~~\forall v\in V;\label{21} \\
&&
a(\tilde{u}_{h},v_{h})= b(f_{h},v_{h})~~~\forall v_{h}\in V^{h}. \label{22}
\end{eqnarray}
Then, for $p_{h}\in W_{f_h}^h$ satisfying $\text{div } p_h = \pi_h \tilde{u}_{h}$,  we have the following revised hypercircle equation:
\begin{eqnarray*}
&&\|\nabla \tilde{u}_{h}-p_{h}\|_{L^2}^2\\
&=&\|\tilde{u}-\tilde{u}_{h}\|_{H^1(\Omega)}^2+\|\nabla\tilde{u}-p_{h}\|_{L^2}^2+\|\tilde{u}-\tilde{u}_{h}\|_{L^2}^2+2((\pi_{h}-I)(\tilde{u}-\tilde{u}_{h}),(\pi_{h}-I)\tilde{u}_{h}).
\end{eqnarray*}
where $I$ is the identity operator.
\end{The}
\proof
Rewriting $\nabla \tilde{u}_{h}-p_{h}$ by $(\nabla \tilde{u}_{h}-\nabla \tilde{u})+(\nabla \tilde{u}-p_{h})$, we have
\begin{eqnarray*}
\|\nabla \tilde{u}_{h}-p_{h}\|_{L^2}^2=\|\nabla \tilde{u}_{h}-\nabla \tilde{u}\|_{L^2}^2+
\|\nabla \tilde{u}-p_{h}\|_{L^2}^2+2(\nabla \tilde{u}_{h}-\nabla \tilde{u},\nabla \tilde{u}-p_{h}).
\end{eqnarray*}
Notice that
\begin{eqnarray*}
&&(\nabla \tilde{u}_{h}-\nabla \tilde{u},\nabla \tilde{u}-p_{h}) =  ( \tilde{u}_{h}-\tilde{u}, -\tilde{u}+\pi_{h}(\tilde{u}_{h}))\\
&=& ( \tilde{u}_{h}-\tilde{u},-\tilde{u}+\tilde{u}_{h}-\tilde{u}_{h}+\pi_{h}(\tilde{u}_{h})) =
\|\tilde{u}_{h}-\tilde{u}\|_{L^2}^2+(\tilde{u}_{h}-\tilde{u},-\tilde{u}_{h}+\pi_{h}(\tilde{u}_{h})).
\end{eqnarray*}
Thus, from the definition of $\pi_{h}$ we get the conclusion.
\endproof


The following theorem gives computable error estimate for $f_h\in  X_{\Gamma}^{h}$.

\begin{The}\label{501}
 Given $f_{h}\in X_{\Gamma}^{h}$, 
let $\tilde{u}\in
V$ and $\tilde{u}_{h} \in V^{h}$ be solutions { to} (\ref{21}) and (\ref{22})
respectively. Then, the following computable error estimate holds:
\begin{eqnarray*}
\|\tilde{u}-\tilde{u}_{h}\|_{H^{1}(\Omega)}&\leq& \kappa_{h}\|f_{h}\|_{b}.
\end{eqnarray*}
Here, $\kappa_{h}$ is defined by
\begin{eqnarray*}
\kappa_{h}:=\max_{f_{h}\in X_{\Gamma}^{h}\setminus \{0\}}\frac{Y(f_{h},p_h, \beta)}{\|f_{h}\|_{b}}
\end{eqnarray*}
where
\begin{eqnarray*}
Y^2(f_{h},p_h,\beta):= (2+\beta+1/\beta)(C_{0}h)^4\|\nabla \tilde{u}_{h} \|_0^2+(1+1/\beta)\|\nabla\tilde{u}_{h}-p_{h}\|_{0}^2 \quad \forall{\beta>0}
\end{eqnarray*}
{and $p_{h}\in W_{f_h}^h$ satisfies $\text{div } p_h = \pi_h \tilde{u}_{h}$}.
\end{The}
\proof
From the hypercircle equation and (\ref{201}), we get
\begin{eqnarray}\label{001}
\|\tilde{u}-\tilde{u}_{h}\|_{H^{1}(\Omega)}^2&\leq& \|\nabla\tilde{u}_{h}-p_{h}\|_{L^2}^2-2((\pi_{h}-I)(\tilde{u}-\tilde{u}_{h}),(\pi_{h}-I)\tilde{u}_{h})\nonumber\\
&\leq& \|\nabla\tilde{u}_{h}-p_{h}\|_{0}^2+2C_{0}h\|\nabla(\tilde{u}-\tilde{u}_{h})\|_{0}\cdot\|(I-\pi_{h})\tilde{u}_{h}\|_{0}\nonumber\\
&\leq& \|\nabla\tilde{u}_{h}-p_{h}\|_{0}^2+2(C_{0}h)^2\|\nabla(\tilde{u}-\tilde{u}_{h})\|_{0}\|\nabla \tilde{u}_{h} \|_0
\end{eqnarray}
Define
$x:= \|\tilde{u}-\tilde{u}_{h}\|_{H^{1}(\Omega)},~
A:= 2( C_{0}h)^2\|\nabla \tilde{u}_{h} \|_0,~
B:=\|\nabla\tilde{u}_{h}-p_{h}\|_{0}$.
By solving the inequality $x^2\leq B^2+Ax$, one can easily deduce that
\begin{equation}
\label{eq:a-posteriori-est}
\|\tilde{u}-\tilde{u}_{h}\|_{H^{1}(\Omega)} \le Y(f_h, p_h, \beta)
\end{equation}
for any $\beta>0$. By further varying $f_h$ in $X_{\Gamma}^h$, we draw the conclusion about $\kappa_{h}$.
\endproof

\begin{remark} The selection of $p_h$ in Theorem \ref{501} is not unique. A proper $p_h$ will be determined in \S\ref{sec:kappa_h_computing}.
In practical computation, since the first term in $Y(f_h, p_h, \beta)$ has higher order convergence, we can take
$\beta > 1$ to have a smaller value of $\kappa_h$.
\end{remark}
%

\section{Explicit { A} Priori Error {Estimates}}\label{sec4}
\subsection {Trace Theorem}
This section is devoted to provide the explicit bound for the constant in the trace theorem.

Let us follow the method in \cite{AVE-2014} to show the explicit value of constants related to trace thereom.

\begin{The}\label{502}
Let $e$ be an edge of triangle element $K$.
Given $u\in V_e(K)$,
we have the following trace theorem
\begin{eqnarray*}
  \|u\|_{L^{2}(e)}&\leq&  0.574 \sqrt{ \frac{|e|}{|K|}}h_K|u|_{H^{1}(K)}\:,
\end{eqnarray*}
where  $V_{e}(K):=\{v\in H^{1}(K) ~~|~~ \int_{e}vds=0 \}$.
\end{The}

\proof
Suppose $\mathbf{P_1}, \mathbf{P_2}, \mathbf{P_3}$ to be the vertices of $K$ and $e:=\mathbf{P_1}\mathbf{P_2}$.
For any $u \in H^{1}(K)$, the Green theorem leads to
$$
\int_K ((x,y)- \mathbf{P_3} ) \cdot \nabla (u^2) \mbox{d}K =
\int_{\partial K} ((x,y)- \mathbf{P_3} )\cdot {\bm{n}} u^2 \mbox{d}s - \int_K 2 u^2 \mbox{d} K\:.
$$
For the term $((x,y)- \mathbf{P_3} )\cdot {\bm{n}}$, we have
\begin{equation}
((x,y)- \mathbf{P_3} )\cdot {\bm{n}} =
\quad
\left\{
\begin{array}{ll}
0, & \mbox{ on }  \mathbf{P_1}\mathbf{P_3}, ~ \mathbf{P_2}\mathbf{P_3}\:, \\
2 {|K|}/{|e|} & \mbox{ on } e\:.
\end{array}
\right.
\end{equation}
Thus,
\begin{align*}
2\frac{|K|}{|e|}\int_e u^2 ds & = \int_K 2 u^2 \mbox{d} K + \int_K ((x,y)- \mathbf{P_3} ) \cdot \nabla (u^2) \mbox{d} K\\
& \le  \int_K 2 u^2 \mbox{d} K + 2h_K \int_K |u| |\nabla u| \mbox{d} K \\
& \le 2\|u\|_{0,K}^2 + 2h_K \|u\|_{0,K} \|\nabla u\|_{0,K}
\end{align*}

Since $u\in V_e(K)$, we have
$$
\int_e u^2  ds  \le \int_e (u-\pi_{h} u)^2  ds  \le
\frac{|e|}{|K|}  \left(\|u-\pi_{h} u\|_{0,K}^2 + h_K \|u-\pi_{h} u\|_{0,K} \|\nabla u\|_{0,K} \right)
$$
By further applying the estimation of $\pi_h$ in (\ref{201}), we have

$$
 \|u\|_{L^{2}(e)}  \le \sqrt{ 1/3.8317^2 +1/3.8317 }\sqrt{ \frac{|e|}{|K|}} ~ h_K ~ \|\nabla u\|_{0,K}
 \le  0.574 \sqrt{ \frac{|e|}{|K|}} ~ h_K ~ \|\nabla u\|_{0,K}
$$

\endproof

\begin{remark}
{
Almost the same result is shown in \cite{AVE-2014}, where general $n$ dimensional element is considered.
Since a sharper bound for $\pi_h$ is utilized here, the constant $0.574$ obtained in Theorem \ref{502} is smaller than the
one in  \cite{AVE-2014} (about 0.648). }
\end{remark}
\begin{remark}
Numerial computation{s indicate} that when the lengths of two edges $P_1P_3$, $P_2P_3$ are fixed as $h$, the constant $C$ in the
 {estimate} $\|u\|_e \le Ch \|\nabla u\|_{0,K}$ $\forall u \in V_e(K)$ will tend to $0$ when the length of the third edge $e:=P_1P_2$ tends to $0$.
However, this behavior of the constant $C$ cannot be deduced from Theorem \ref{502}.
\end{remark}

\subsection{{Explicit A} Priori Error {Estimates}}


\begin{The}\label{503}
Let $u$ and $\tilde{u}$
be solutions { to} (\ref{11}) and (\ref{21}){,} respectively, {with $f_h$} taken as $f_h:=\pi_{h,\Gamma}f$. Then, the following error estimate holds:
\begin{eqnarray*}\label{31}
\|u-\tilde{u}\|_{H^{1}(\Omega)}
\leq C_{1}(h)\|(I-\pi_{h,\Gamma})f\|_{b}.
\end{eqnarray*}
where $C_{1}(h)=\max\limits_{e\in E_{h,\Gamma}}\{  {0.574}\sqrt{ \frac{|e|}{|K|}}h_K\}$.
\end{The}
\proof
Setting $v=u-\tilde{u}$ in (\ref{11}) and (\ref{21}), we have
\begin{eqnarray*}
a(u-\tilde{u}, u-\tilde{u})&=&b(f-f_{h},u-\tilde{u})\\
&=&b((I-\pi_{h,\Gamma})f,(I-\pi_{h,\Gamma})(u-\tilde{u})).
\end{eqnarray*}
From the Schwartz inequality and Theorem \ref{502}, we get
\begin{eqnarray*}
\|u-\tilde{u}\|_{H^{1}(\Omega)}^2&\leq&\|(I-\pi_{h,\Gamma})f\|_{b}\|(I-\pi_{h,\Gamma})(u-\tilde{u})\|_{b}\\
&\leq&C_{1}(h)\|(I-\pi_{h,\Gamma})f\|_{b}|u-\tilde{u}|_{H^{1}(\Omega)}
\end{eqnarray*}
which {implies} the conclusion.
\endproof


%
%

Now, we are ready to {formulate and prove} the explicit {a} priori error {estimate}.

\begin{The}
Let $u$ and $u_{h}$
be solutions {to} (\ref{11}) and (\ref{12}){,} respectively. Then, the following error estimates hold:
\begin{eqnarray}
\|u-u_{h}\|_{H^{1}(\Omega)}&\leq& M_{h}\|f\|_{b}\label{50}\\
\|u-u_{h}\|_{b}&\leq& M_{h}^2\|f\|_{b}\label{51}
\end{eqnarray}
with $M_{h}:=\sqrt{(C_{1}(h))^2+\kappa_{h}^2}.$
\end{The}
\proof
The estimation in (\ref{50}) can be obtained by applying Theorems \ref{501} and \ref{503},
\begin{eqnarray*}
\|u-u_{h}\|_{H^{1}(\Omega)}&\leq&   \|u-\tilde{u}\|_{H^{1}(\Omega)}+\|\tilde{u}-\tilde{u}_{h}\|_{H^{1}(\Omega)}\\
&\leq& C_{1}(h)\|(I-\pi_{h,\Gamma})f\|_{b}+\kappa_{h}\|f_{h}\|_{b}\\
&\leq& \sqrt{(C_{1}(h))^2+\kappa_{h}^2}\|f\|_{b}\:.
\end{eqnarray*}

The error estimate (\ref{51}) can be obtained by the Aubin--Nitsche duality technique.
\endproof

{\begin{remark}
The result (\ref{50}) of Theorem 4.3 provides an explicit {\it a priori} error estimation for the FEM solutions,
which is based on the {\it  a posteirori} error estimation in (\ref{eq:a-posteriori-est}).
Notice that in (\ref{eq:a-posteriori-est}), by taking any explicit $p_h$ and $\beta$, we have the following explicit {\it  a posteriori} bound for the FEM solution.
\begin{equation}
\label{eq:full-a-post-error}
\|u-{u}_{h}\|_{H^{1}(\Omega)} \leq C_{1}(h)\|(I-\pi_{h,\Gamma})f\|_{b} + Y(f_h, p_h, \beta).
\end{equation}
Similar results about {\it a posteriori } error estimation can be found in \cite{AVE,AVE-2014, Kikuchi-Saito}:
In \cite{ AVE, Kikuchi-Saito}, the homogeneous Dirichlet boundary condition is considered;
In \cite{AVE-2014}, the nonhomogeneous Neumann boundary condition is considered and (\ref{eq:full-a-post-error})
can be regarded as a special case of \cite{AVE-2014}.
\end{remark}
}

\subsection{Computation of $\kappa_{h}$ }\label{sec:kappa_h_computing} The quantity $\kappa_{h}$ is evaluated in two steps. \\

First, for fixed $f_{h}$, we deduce explicit forms of $\tilde{u}_{h}\in V^{h}$ and
$p_{h}\in W^{h}$, which appear in the definition of $Y(f_{h},p_h, \beta)$.
 According to the standard theories of the conforming FEM and the Raviart-Thomas FEM; see, e.g., \cite{BrezziFortin},
we solve the following two problems:

\begin{itemize}
\item [(a)] Find $\tilde{u}_{h}\in V^{h}$ such that
\begin{eqnarray*}\label{59}
a(\tilde{u}_{h}, v_{h})=b(f_{h},v_{h}) \quad \forall v_{h}\in V^{h}.
\end{eqnarray*}
\item [(b)] Let $\tilde{u}_{h}$ be the solution of (a). Find $p_{h}\in W^{h}_{f_h}$ and $\rho_{h}\in X^{h}$, $c\in \mathbf{R}$ such that
\begin{equation*}\label{020}
\left\{
\begin{array}{rcl}
(p_{h},\tilde{p}_{h})+(\rho_{h},\mbox{div } \tilde{p}_{h}) + (\rho_h, d) &=& { 0 } \ \ \ \ \  \forall \tilde{p}_{h}\in W_0^{h}, ~ {\forall d \in \mathbf{R}}{,} \\
(\mbox{div } p_{h},\tilde{q}_{h}) + ({c}, \tilde{q}_{h}) &=& (\pi_{h}(\tilde{u}_{h}), \tilde{q}_{h})\ \ \ \ \  \forall \tilde{q}_{h}\in X^{h}{,}
\end{array}
\right.
\end{equation*}
where ${W_0^{h}:=\{p_{h}\in W^{h} \: | \: p_{h}\cdot {\bm{n}}= 0 \in X_{\Gamma}^{h}~\text{on}~\Gamma \}}\:.$

\end{itemize}

\vskip 0.5cm
Notice that the solution $p_h$ of (b) is depending on $f_h$. Let us rewrite $Y(f_h,p_h, \beta)$ as $Y(f_h,\beta)$.
Second, we find $f_{h}$ that maximizes the
value of $Y(f_h,\beta)/\|f_{h}\|_{b}$ by solving an eigenvalue problem.
By using the solutions of (a) and (b), $Y(f_h,\beta)$ {and} $\|f_h\|{_{b}}$ can be formulated by
$$
Y^2 (f_h,\beta)= x^TAx {~~\mbox{and}} ~~ \|f_h\|{_{b}}^2 = x^TBx \:.
$$
where $x$ is the coefficent vector of $f_h$ with respect to the basis of $X^h_{\Gamma}$ and
$A$, $B$ are symmetric matrices to be determined upon the selection of basis of FEM spaces.
Thus, the value of $\kappa_h^2$ is given by the maximum eigenvalue of the problem
$$
Ax=\lambda {B}x\:.
$$
For detailed {solution} of {this} eigenvalue problem, {we} refer to \cite{LO}{, where an} analogous  {problem is described}.

\section{Numerical Examples}\label{sec5}

In this section, {several
numerical tests are presented. The constant $\kappa_{h}$ is computed for problem (\ref{1})
and four different domains. For each domain a sequence of uniformly refined
finite element meshes is considered. If $\kappa_{2h}$ and $\kappa_{h}$ are computed on two consecutive meshes then the convergence rate is estimated numerically as$$
\kappa_h\mbox{-rate}:=\log(\kappa_{2h}/\kappa_{h})/\log2\:.
$$}
\subsection{{T}he unit square }
We consider the problem (\ref{1}) on the unit square domain $\Omega=(0,1)\times (0,1)$.
In the numerical experiment, we set $\beta=0.1, 1, 10, 100$ and $1000$.
The dependency of $\kappa_h$ on $\beta$ is displayed in Figure 1,
which {illustrates} that {larger $\beta$} gives smaller $\kappa_h$. {However,} the definition of
$Y(f_h,\beta)$ {clearly shows that $\beta$ cannot be too large}.

{Computed} quantities $\kappa_{h}$, $C_{1}(h)$, and $M_{h}$ {for the case} $\beta=100$ are shown in Table 1.
{ The estimated} convergence {rate} of $\kappa_h$, denoted by $\kappa_{h}$-{rate}, is close to $0.5$.

\begin{table}[h]
\begin{center}
\caption{{ Computed q}uantities for the square { and} $\beta=100$}
\begin{tabular}{|c|c|c|c|c|}\hline
$h$&$\kappa_{h}$&$C_{1}(h)$&$M_{h}$&$\kappa_{h}$-{ rate}\\
\hline
${\sqrt{2}}/{4}$&~0.4143&~~ 0.574  &~~0.7079 &~ - \\ \hline
${\sqrt{2}}/{8}$&~0.2973&~~0.4059 &~~0.5031 &~0.4788 \\\hline
${\sqrt{2}}/{16}$&~0.2110&~~0.2870&~~0.3562&~0.4947\\\hline
${\sqrt{2}}/{32}$&~0.1493&~~ 0.2029&~~0.2519&~0.4990 \\\hline
\end{tabular}
\end{center}
\end{table}

\begin{figure*}
\begin{center}
  \includegraphics[width=0.7\textwidth]{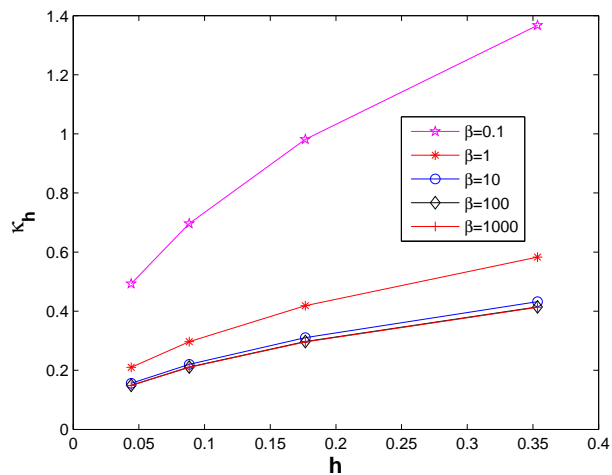}
\caption{The dependence of $\kappa_h$ on $\beta$ (unit square)}
\label{fig:1}       
\end{center}
\end{figure*}

\subsection{{ Right} triangle{, equilateral triangle,} and the L-shape domain}

{ In this example, three} domains are {considered, namely, the isosceles} right triangle {with unit legs}, the unit { equilateral} triangle, and the L-shaped domain
$\Omega=(0, 1)\times (0, 1)\setminus [1/2, 1]\times [1/2, 1]$.
The results {for $\beta$=100 are displayed in Tables 2-4, respectively. For all domains the convergence rate of $\kappa_{h}$} is close to $0.5$.

\begin{table}[h]
\begin{center}
\caption{{Computed q}uantities for the {isosceles} right triangle { and} $\beta=100$}
\begin{tabular}{|c|c|c|c|c|}\hline
$h$&$\kappa_{h}$&$C_{1}(h)$&$M_{h}$&$\kappa_{h}$-{ rate}\\
\hline
${\sqrt{2}}/{4}$&~0.4448&~~0.6826&~~0.8147 &~ - \\\hline
${\sqrt{2}}/{8}$&~0.3107&~~0.4827&~~0.5741 &~0.5176\\\hline
${\sqrt{2}}/{16}$&~0.2197&~~ 0.3413 &~~0.4059&~0.5000\\\hline
${\sqrt{2}}/{32}$&~0.1554&~~ 0.2413&~~0.2870&~0.4995\\\hline
\end{tabular}
\end{center}
\end{table}

\begin{table}[h]
\begin{center}
\caption{ { Computed q}uantities for the { equilateral} triangle { and} $\beta=100$}
\begin{tabular}{|c|c|c|c|c|}\hline
$h$&$\kappa_{h}$&$C_{1}(h)$&$M_{h}$&$\kappa_{h}$-{ rate}\\
\hline
$1/4$&~0.3783&~~0.4361&~~0.5773&~ - \\\hline
$1/8 $&~0.2696&~~0.3084&~~0.4096&~0.4887\\\hline
$1/16$&~0.1909&~~0.2181&~~ 0.2898&~0.4980 \\\hline
$1/32$&~0.1350&~~0.1542&~~0.2049&~0.4999\\\hline
\end{tabular}
\end{center}
\end{table}

\begin{table}[h]
\begin{center}
\caption{{ Computed q}uantities for the L-shape domain { and} $\beta=100$}
\end{center}
\begin{center}
\begin{tabular}{|c|c|c|c|c|}\hline
$h$&$\kappa_{h}$&$C_{1}(h)$&$M_{h}$&$\kappa_{h}$-{ rate}\\
\hline
${\sqrt{2}}/{4}$&~0.4872&~~0.574  &~~0.7529&~ - \\\hline
${\sqrt{2}}/{8}$&~0.3432&~~0.4059&~~0.5315&~0.5055\\\hline
${\sqrt{2}}/{16}$&~0.2439&~~0.2870&~~ 0.3766&~0.4928\\\hline
${\sqrt{2}}/{32}$&~0.1734&~~0.2029&~~0.2669&~0.4922\\\hline
\end{tabular}
\end{center}
\end{table}

\section{\bf Conclusion}

In this paper, by applying the technique of the hypercircle equation, we successfully construct the explicit
a priori error {estimate} for the FEM solution of { nonhomogeneous Neumann problems}. By following the framework proposed by the second author
in \cite{Liu}, the a priori error {estimate} obtained here can be used in bounding eigenvalues of the Steklov type eigenvalue problems. {The expected rate of convergence of
$M_h$ is $1$ in case the solution is smooth enough. In this paper, only the $H^1$ regularity is required in the analysis,
and both the theoretical results, see Theorem 4.2, and numerical tests
confirm the suboptimal convergence rate $0.5$ for $M_h$ as well as $\kappa_h$. It is an interesting problem that
whether the rate of convergence can be improved or not, for general $f\in L_2(\partial \Omega)$.}

\vskip 0.5cm

{\small
\indent

{\em Authors' addresses}:\\

{\em Qin LI}, School of Science, Beijing Technology and Business University, Beijing 100048, P. R. China,
 e-mail: \texttt{liqin@\allowbreak lsec.cc.ac.cn}; \\

 {\em Xuefeng LIU} (correspoding author), Graduate School of Science and Technology, Niigata University, 8050 Ikarashi 2-no-cho, Nishi-ku, Niigata City, Niigata 950-2181 Japan, e-mail: \texttt{xfliu@\allowbreak math.sc.niigata-u.ac.jp}.

}

\end{document}